\documentstyle[12pt]{article}

\def\re{\par\hang\textindent}

\def\be{\begin{equation}}
\def\ee{\end{equation}}

\def\qed{$\hfill\rule[0cm]{1.5mm}{3mm}$}
\begin{document}

\vspace*{1cm}

\centerline{\large\bf Difference independence of the Riemann zeta
function} \vskip 0.6cm \centerline{Yik-Man Chiang$^{\ddag *}$
\footnote{*  This research was supported in part by the Research
Grant Council of the Hong Kong Administrative Region, China
(HKUST6135/01P)}\ and \ Shao-Ji Feng$^{\dag \ddag **}$\footnote{**
The work of this author is supported by the National Natural Science
Foundation of China (Grant No. 10501044), in part by the Research
Grant Council of the Hong Kong Administrative Region, China
(HKUST6135/01P) and HKUST PDF Matching Fund}} \vskip 3mm
\centerline{$\dag$ Academy of Mathematics and Systems Science,}
\centerline{ Chinese Academy of Sciences, Beijing, 100080,\ P. R.
China.} \centerline{e-mail: fsj@amss.ac.cn} \vskip 3mm
\centerline{$\ddag$ Department of Mathematics, Hong Kong University
of} \centerline{Science and Technology, Clear Water Bay, Kowloon,}
\centerline{Hong Kong,\ P. R. China.} \centerline{e-mail:
machiang@ust.hk} \vskip 1cm
\begin{quote}
\small {\bf Abstract.\ } It is proved that the Riemann zeta function
does not satisfy any nontrivial algebraic difference equation whose
coefficients are meromorphic functions $\phi$ with Nevanlinna
characteristic
satisfying $T(r, \phi)=o(r)$ as $r\rightarrow \infty$.\\[3mm]
\\[3mm]
\end{quote}

\vskip 5mm

\def\theequation{1.\arabic{equation}}
\setcounter{equation}{0} \noindent{\large\bf 1. \ Introduction.}
\vskip 8mm

Let $\zeta(s)$ denote the Riemann zeta function, and $\Gamma(s)$
denote the Eular gamma function in the complex plane. These two
important special functions are related by the functional equation
(see [7, 16]): \be \zeta(1-s)=2^{1-s}\pi^{-s}\cos(\frac12\pi
s)\Gamma(s)\zeta(s). \ee

A classical theorem of H\"{o}lder [6] concerning $\Gamma(s)$ states
that $\Gamma(s)$ can not satisfy any algebraic differential equation
whose coefficients are rational functions [i.e. any equation of the
form $f(s, y, y^{\prime}, \cdots, y^{(n)})=0$, where $n$ is a
nonnegative integer and where $f$ is a polynomial in $y, y^{\prime},
\cdots, y^{(n)}$ whose coefficients are rational functions of $s$].
Other proofs of this theorem were given in [3, 9, 11, 12]. Bank and
Kaufman [1] generalized the theorem to coefficients being
meromorphic functions $\phi$ with Nevanlinna characteristic (see
(1.5) below for its definition) satisfying $T(r, \phi)=o(r)$ as
$r\rightarrow \infty$.

The question of the differential independence of $\zeta(s)$ was
touched upon by Hilbert in 1900 (see [5]). He conjectured that
$\zeta(s)$ and other functions of the same type do not satisfy
algebraic differential equations with rational coefficients. It
follows from Hilbert's report that, based on H\"{o}lder's theorem on
the algebraic differential independence of $\Gamma(s)$, he could
prove the algebraic differential independence of $\zeta(s)$.
Hilbert's conjectures were proved in [10, 13], see also [14, 15].

It is well known that the gamma function satisfies the following
difference equation \be \Gamma(s+1)=s\Gamma(s).\ee A natural
question is to ask whether $\zeta(s)$ satisfies any algebraic
difference equation, or doing $\zeta(s)$ difference independent in
short? More precisely, if $\zeta(s+s_0), \zeta(s+s_1), \cdots,
\zeta(s+s_m)$ satisfy any algebraic equation or not? here $m$ is a
nonnegative integer, $s_i,\ i=0, 1, \cdots, m$ are distinct {\bf
complex} numbers. The special case that $s_0, s_1,\cdots, s_m$ are
{\bf real} numbers was studied by Ostrowski [13]. He proved that
$\zeta(s+s_0), \zeta(s+s_1), \cdots, \zeta(s+s_m)$ can not satisfy
any algebraic equation with rational coefficients. For `small'
$s_i$'s, Voronin [17] proved the following.

\vskip 3mm

{\bf Theorem A.} \ Let $m$ be a nonnegative integer, $s_i\in {\bf
C},\ i=0, 1, \cdots, m$, $s_i\not=s_j, 0\leq i<j\leq m$, if \be
|s_i|<\frac14,\ \ \ i=0, 1, \cdots, m,\ee and \be f(\zeta(s+s_0),
\zeta(s+s_1), \cdots, \zeta(s+s_m))=0\ee identically in $s\in C$,
where $f(z_0, z_1, \cdots, z_m)$ is a continuous function, then $f$
is identically zero.

\vskip 3mm

 Theorem A is not always true again for
`large' $s_i$'s, since the vectors $(\zeta(s),\zeta(s+s_0)$ are not
dense in ${\bf C}^2$ for any given complex $s_0$ with $\Re s_0>1$,
as follows from elementary properties of the zeta and gamma
functions.

   In order to state our main result, we recall the standard notations of the Nevanlinna theory. Let
$\phi$ be a meromorphic function on the complex plane, the
Nevanlinna characteristic $T(r,\phi)$ of $\phi$ for $r\geq 0$ is
defined by \be T(r,\phi)=N(r,\phi)+m(r,\phi),\ee where \be
N(r,\phi)=\int_0^r\frac{n(t,\phi)-n(0,\phi)}{t}dt+n(0,\phi)\log r\ee
is the pole counting function, $n(r,\phi)$ is the number of poles
(counting multiplicities) of $\phi$ in the disc $\{z:\ |z|\leq r\}$,
the proximity function $m(r,\phi)$ is given by \be
m(r,\phi)=\frac{1}{2\pi}\int_0^{2\pi}\log^+|\phi(re^{{\mathrm
i}\theta})|d\theta,\ee where $\log^+x=\max\{0,\log x\}$. We have

 \vskip 3mm

{\bf Theorem 1.} \  The Riemann zeta function $\zeta(s)$ does not
satisfy any algebraic difference equation whose coefficients are
meromorphic functions $\phi$ with Nevanlinna characteristic
satisfying $T(r, \phi)=o(r)$ as $r\rightarrow \infty$. That is, if
\be f(s, \zeta(s+s_0), \zeta(s+s_1), \cdots, \zeta(s+s_m))=0\ee
holds for all $s\in {\bf C}$, where $m$ is a nonnegative integer,
$s_i\in C,\ i=0, 1, \cdots, m$, are distinct complex numbers, and
$f$ is a polynomial in $\zeta(s+s_0), \zeta(s+s_1), \cdots,
\zeta(s+s_m)$ whose coefficients are meromorphic functions $\phi(s)$
with Nevanlinna characteristic satisfying $T(r, \phi)=o(r)$ as
$r\rightarrow \infty$, then $f$ is identically zero.

\vskip 3mm

It is well known that a meromorphic function $\phi$ is rational if
and only if $T(r,\phi)=O(\log r)$ as $r\rightarrow\infty$. Theorem 1
implies that $\zeta(s)$ can not satisfy any algebraic difference
equation whose coefficients are rational functions.

We recall that the Nevanlinna order of meromorphic function $\phi$
is defined by \be \rho(\phi)=\limsup_{r\rightarrow\infty}\frac{\log
T(r,\phi)}{\log r}.\ee It follows form Theorem 1 that

\vskip 3mm

{\bf Corollary 1.} \ The Riemann zeta function $\zeta(s)$ does not
satisfy any algebraic difference equation whose coefficients are
meromorphic functions $\phi$ with Nevanlinna order $\rho(\phi)<1$.

\vskip 3mm

It is easily seen that the Nevanlinna order of $\zeta(s)$ is $
\rho(\zeta(s))=1$. By a similar calculation, we have \be
\rho(\frac{\zeta(s+1)}{\zeta(s)})=1.\ee It follows that Corollary 1
is best possible in the sense that the condition $\rho(\phi)<1$ can
not be relaxed to $\rho(\phi)\leq 1$.

 \vskip 8mm

\def\theequation{2.\arabic{equation}}
\setcounter{equation}{0} \noindent {\large\bf 2.\ Proof of the
constant coefficients case of Theorem 1.} \vskip 8mm

In this section, we prove the constant coefficients case of Theorem
1.

\vskip 3mm

{\bf Theorem 2.} \  The Riemann zeta function $\zeta(s)$ does not
satisfy any algebraic difference equation with constant
coefficients. That is, if \be f(\zeta(s+s_0), \zeta(s+s_1), \cdots,
\zeta(s+s_m))=0\ee for all $s\in {\bf C}$, where $m$ is a
nonnegative integer, $s_i\in {\bf C},\ i=0, 1, \cdots, m$, are
distinct complex numbers, $f$ is a polynomial in $\zeta(s+s_0),
\zeta(s+s_1), \cdots, \zeta(s+s_m)$ whose coefficients are complex
constants, then $f$ is identically zero.

\vskip 3mm

To prove Theorem 2, we need the following lemma.

\vskip 3mm

{\bf Lemma 1.} \ Let $b_i,\ s_i\in C,\ i=0, 1, \cdots, m$,
$s_i\not=s_j, 0\leq i<j\leq m$. Suppose that there exists an integer
$p_0>0$ such that for all prime numbers $p>p_0$, \be
\sum_{i=0}^mb_ip^{-s_i}=0.\ee Then \be b_i=0, \ \ \ i=0,\ 1,\
\cdots,\ m.\ee

{\bf Proof.} We construct an entire function of exponential type \be
g(z)=\sum_{i=0}^mb_ie^{-s_iz},\ee It is therefore not difficult to
see that $g(z)$ can not have more that $cr$ zeros in a disc of
radius $r$, where $c$ is a constant, unless it vanishes identically.
It follows from the hypothesis (2.2) that $g(z)$ vanishes at $z=\log
p$ for prime $p>p_0$. We deduce from the Prime Number Theorem that
the number of zeros of $g(z)$ in a disc of radius $r$ is at least
$$ \frac{e^r}{r}+o(\frac{e^r}{r})\geq cr$$ as $r\rightarrow+\infty$.
Therefore $g(z)$ is identically zero and then (2.3) follows.

\vskip 3mm

{\bf Proof of Theorem 2.} \ Assume that $\zeta(s)$ satisfy \be
\sum_{i=1}^NP_i\zeta(s+s_0)^{k_0(i)}\zeta(s+s_1)^{k_1(i)}\cdots\zeta(s+s_m)^{k_m(i)}
=0,\ee where $ P_i\ ( i=1,\ 2,\ \cdots, \ N)$ are complex constants,
not all are zero, $K(i)=(k_0(i), k_1(i), \cdots, k_m(i)), \ i=1,\
2,\ \cdots, \ N$, are multi-indices with all indices being
nonnegative integers, such that \be K(i)\not=K(j), \ \ \ \ 1\leq
i<j\leq N. \ee Here $K(i)=K(j)$ means $ k_l(i)=k_l(j)$ for each $l,
0\leq l\leq m$. Let \be
\sum_{n=1}^{\infty}\frac{A_i(n)}{n^s}:=\zeta(s+s_0)^{k_0(i)}\zeta(s+s_1)^{k_1(i)}\cdots\zeta(s+s_m)^{k_m(i)},\
\ \ i=1,\ 2,\ \cdots, \ N,\ee where the Dirichlet series is
convergent in the region $1-\min\{\Re s_0, \Re s_1, \cdots, \Re
s_m\}<\Re s<\infty$, then by the Uniqueness Theorem for Dirichlet
series, we have \be \sum_{i=1}^NP_iA_i(n)=0,\ \ \ \ n=1,\ 2,\ 3,
\cdots.\ee

If there exists an $i^{\prime}$,  $1\leq i^{\prime}\leq N $ such
that $K(i^{\prime})=(0, 0, \cdots, 0)$, then there exists an $i$,
$i\not=i^{\prime}$ such that $P_i\not=0$. (Otherwise we would have
$P_{i^{\prime}}\not=0$ and $P_{i^{\prime}}A_{i^{\prime}}(1)=0$. But
$A_{i^{\prime}}(1)=1$, so we have a contradiction.) Since
$A_i(1)=1$, so \be A_{i^{\prime}}(n)=0,\ \ \ \ n=2,\ 3, \ 4,
\cdots,\ee we have \be \sum_{1\leq i\leq N, i\not=i^{\prime}}
P_iA_i(n)=0,\ \ \ \ n=2,\ 3,\ 4, \cdots\ee with $P_i,\ 1\leq i\leq
N,\ i\not=i^{\prime}$ not all zero.

So, we may assume without loss of generality that \be K(i)\not=(0,
0, \cdots, 0), \ \ \ \ 1\leq i\leq N,\ee and \be
\sum_{i=1}^NP_iA_i(n)=0,\ \ \ \ n=2,\ 3,\ 4, \cdots\ee hold with
$P_i,\ i=1, 2, \cdots, N$, not all equal to zero.

Let $L$ be an arbitrary positive integer, $p_1, p_2, \cdots, p_L$ be
distinct prime numbers. By (2.7), we have for $i=1, 2, \cdots, N$,
\be A_i(p_1p_2\cdots p_L)=\prod_{1\leq l\leq
L}(k_0(i)p_l^{-s_0}+k_1(i)p_l^{-s_1}+\cdots+k_m(i)p_l^{-s_m}).\ee It
follows from (2.8) and (2.13) that \be \sum_{i=1}^NP_i\prod_{1\leq
l\leq
L}(k_0(i)p_l^{-s_0}+k_1(i)p_l^{-s_1}+\cdots+k_m(i)p_l^{-s_m})=0.\ee

Let $p_1, p_2, \cdots, p_{L-1}$ be fixed, and $p_L$ varies all over
$\{p\ | p\ \mbox{is a prime number},\ p>\max_{1\leq l\leq
L-1}p_l\}$. We then have for all prime numbers $ p>\max_{1\leq l\leq
L-1}p_l$, \be \sum_{q_1=0}^m\Big[\sum_{i=1}^NP_i\prod_{1\leq l\leq
L-1}(k_0(i)p_l^{-s_0}+k_1(i)p_l^{-s_1}+\cdots+k_m(i)p_l^{-s_m})
k_{q_1}(i)\Big]p^{-s_{q_1}}=0.\ee Thus by Lemma 1, we have for
$q_1=0, 1, \cdots, m$, \be \sum_{i=1}^NP_i\prod_{1\leq l\leq
L-1}(k_0(i)p_l^{-s_0}+k_1(i)p_l^{-s_1}+\cdots+k_m(i)p_l^{-s_m})
k_{q_1}(i)=0.\ee

Let $p_1, p_2, \cdots, p_{L-2}$ be fixed, and $p_{L-1}$ varies all
over $\{p\ | p\ \mbox{is a prime number},\ p>\max_{1\leq l\leq
L-2}p_l\}$. We then have for $q_1=0, 1,\cdots, m$, for all prime
numbers $p>\max_{1\leq l\leq L-2}p_l$, \be
\sum_{q_2=0}^m\Big[\sum_{i=1}^NP_i\prod_{1\leq l\leq
L-2}(k_0(i)p_l^{-s_0}+k_1(i)p_l^{-s_1}+\cdots+k_m(i)p_l^{-s_m})k_{q_1}(i)
k_{q_2}(i)\Big]p^{-s_{q_2}}=0.\ee Thus by Lemma 1 again, we have for
$q_1, q_2=0, 1, \cdots, m$, \be \sum_{i=1}^NP_i\prod_{1\leq l\leq
L-2}(k_0(i)p_l^{-s_0}+k_1(i)p_l^{-s_1}+\cdots+k_m(i)p_l^{-s_m})
k_{q_1}(i)k_{q_2}(i)=0.\ee

Continuing this procedure, we finally have for $L=1, 2, 3, \cdots$,
\be \sum_{i=1}^NP_ik_{q_1}(i)k_{q_2}(i)\cdots k_{q_L}(i)=0,\ \ \
q_1, q_2, \cdots, q_L=0, 1, \cdots, m.\ee Then
\begin{eqnarray}
 &&\sum_{i=1}^NP_i(k_0(i)+k_1(i)x+\cdots+k_m(i)x^m)^L\nonumber\\
&=&\sum_{i=1}^NP_i\sum_{0\leq q_1, q_1, \cdots, q_L\leq m}k_{q_1}(i)k_{q_2}(i)\cdots k_{q_L}(i)x^{q_1+q_2+\cdots+q_L}\nonumber\\
&=&\sum_{0\leq q_1, q_1, \cdots, q_L\leq m}x^{q_1+q_2+\cdots+q_L}\sum_{i=1}^NP_ik_{q_1}(i)k_{q_2}(i)\cdots k_{q_L}(i)\nonumber\\
&=&0,\ \ \ \ \ \ \ \ \forall x\in \textbf{R},\ \ \ L=1, 2, 3,
\cdots.\end{eqnarray} Denote \be
H_i(x)=k_0(i)+k_1(i)x+\cdots+k_m(i)x^m,\ \ \ \ \ i=1, 2, \cdots N\ee
and \be H_0(x)=0.\ee Notice that $P_i,\ i=1, 2, \cdots, N$, not all
equal zero, by a well-known theorem in determinants due to
Vandermonde, we have $\forall x\in R$, there exists $0\leq
i_1(x)<i_2(x)\leq N$ such that \be H_{i_1(x)}(x)=H_{i_2(x)}(x),\ee
here $i_1(x), i_2(x)$ means the index may depend on $x$. Let \be
E_{ij}=\{x\ | H_i(x)=H_j(x)\}, \ \ \ \ 0\leq i<j\leq N,\ee then \be
\bigcup_{0\leq i<j\leq N}E_{ij}=R.\ee Hence there must exists $0\leq
\hat{i}<\hat{j}\leq N$ such that $E_{\hat{i}\hat{j}}$ is a infinite
set. That is, the polynomial \be H_{\hat{i}}(x)-H_{\hat{j}}(x)\ee
have infinitely many zeros. Then $H_{\hat{i}}(x)-H_{\hat{j}}(x)$
must be the zero-polynomial, which is contradict with (2.7) and
(2.12). The proof is complete. \qed

\vskip 8mm
\def\theequation{3.\arabic{equation}}
\setcounter{equation}{0} \noindent{\large\bf 3.\ Proof of Theorem
1.} \vskip 8mm

We need the following well-known lemma of Cartan [2] (see also [8]).

\vskip 3mm

{\bf Lemma 2.} \ Let $z_1, z_2, \cdots, z_p$ be any finite
collection of complex numbers, and let $B>0$ be any given positive
number. Then there exists a finite collection of closed disks
$D_1, D_2, \cdots, D_q$ with correspending radii $r_1, r_2,
\cdots, r_q$ that satisfy \be r_1+r_2+\cdots+r_q=2B,\ee such that
if $z\notin D_j$ for $j=1,2,\cdots,q$, then there is a permutation
of the points $z_1,z_2,\cdots,z_p$, say,
$\hat{z_1},\hat{z_2},\cdots,\hat{z_p}$, that satisfies \be
|z-\hat{z_{\mu}}|>B\frac{\mu}{p}, \ \ \ \ \ \ \ \ \mu=1, 2,
\cdots, p,\ee where the permutation may depent on $z$.

\vskip 3mm

The following two lemmas may have their independent interest.

\vskip 3mm

{\bf Lemma 3.} \  Let $\phi$ be a meromorphic function, and let
$\alpha>1,\ \varepsilon>0$ be given real constant. Then there
exists a set $E\subset (1,\infty)$ that has {\sl upper logarithmic
density} \be
\delta(E):=\limsup_{x\rightarrow+\infty}\frac{\int_{E\cap
(1,x]}\frac{1}{t}dt}{\log x}<\varepsilon \ee and constant $A>0$
such that for all $z$ satifying $|z|=r\notin [0,1]\cup E$, we have
\be  |\phi(z)|\leq e^{AT(\alpha r,\phi)} .\ee

\vskip 3mm

{\bf Remark.} \ Although An upper bound estimate similar to (3.4)
for the meromorphic function $\phi(z)$ outside an exceptional set
should have been recorded in the literature, the authors is not able
to find such a reference.

\vskip 3mm

{\bf Proof.} Assume $\phi$ is not identically zero. Let $(a_{\nu})_{\nu\in N}$, resp. $(b_{\mu})_{\mu\in N}$, denote the sequence of all zeros, resp. all poles,
of $\phi$, with due account of multiplicity,  and $\Re s$ means the real part of $s$. By the Poisson-Jensen formula, we have for $|z|=r<R$,
\begin{eqnarray}
\log |\phi(z)|&=&\frac{1}{2\pi}\int_0^{2\pi}\log |\phi(Re^{{\mathrm i}\theta})|\Re(\frac{Re^{{\mathrm i}\theta}+z}{Re^{{\mathrm i}\theta}-z})d\theta\nonumber\\
&&-\sum_{|a_{\nu}|<R}\log|\frac{R^2-\bar{a_{\nu}}z}{R(z-a_{\nu})}|+
\sum_{|b_{\mu}|<R}\log|\frac{R^2-\bar{b_{\mu}}z}{R(z-b_{\mu})}|\nonumber\\
&\leq& \frac{R+r}{R-r}\cdot\frac{1}{2\pi}\int_0^{2\pi}|\log
|\phi(Re^{{\mathrm i}\theta})||d\theta+
\sum_{|b_{\mu}|<R}\log|\frac{R^2-\bar{b_{\mu}}z}{R(z-b_{\mu})}|\nonumber\\
&\leq& \frac{R+r}{R-r}(m(R,\phi)+m(R,\frac{1}{\phi}))+
\sum_{|b_{\mu}|<R}\log\frac{R+r}{|z-b_{\mu}|},\end{eqnarray} where
we have used the following estimate (see [4])  \be
|\frac{R^2-\bar{a}z}{R(z-a)}|>1, \ \ \mbox{for}\ |z|<R, \ |a|<R. \ee
Let $R=\alpha^{\frac13} r$, $\alpha>1$, we get
\begin{eqnarray} \log |\phi(z)|&\leq& A_1 T(\alpha^{\frac13} r,\phi)+\sum_{|b_{\mu}|<\alpha^{\frac13} r}\log\frac{(\alpha^{\frac13}+1)r}{|z-b_{\mu}|}\nonumber\\
&\leq& A_1 T(\alpha r,\phi)+\sum_{|b_{\mu}|<\alpha^{\frac13} r}\log\frac{(\alpha^{\frac13}+1)r}{|z-b_{\mu}|},
\end{eqnarray}
where $A_1>0$ is a constant.

Now we estimate $\sum_{|b_{\mu}|<\alpha^{\frac13}
r}\log\frac{(\alpha^{\frac13}+1)r}{|z-b_{\mu}|}$. We suppose that
$h$ is a fixed nonnegative integer, and that $z$ is confined to the
annulus \be \alpha^{\frac{h}{3}}\leq |z|=r\leq
\alpha^{\frac{h+1}{3}}.\ee Set $p=n(\alpha^{\frac{h+2}{3}},\phi)$,
$B=\varepsilon \alpha^{\frac{h}{3}}$, and apply Lemma 2 to the
points $b_1,b_2,\cdots, b_p$, we obtain that there exists a finite
collection of closed disks $D_1, D_2, \cdots, D_q$, whose radii has
a total sum equal to $2B$, such that if $z\notin D_j$ for
$j=1,2,\cdots,q$, then there is a permutation of the points
$b_1,b_2,\cdots,b_p$, say, $\hat{b_1},\hat{b_2},\cdots,\hat{b_p}$,
such that the inequalities \be |z-\hat{b_{\mu}}|>B\frac{\mu}{p}, \ \
\ \ \ \ \ \ \mu=1, 2, \cdots, p\ee hold. Notice that here $q$ and
$D_1, D_2, \cdots, D_q$ depend on $p$ and then depend on $h$. Hence
if $z\notin D_j$ for $j=1,2,\cdots,q$, we have form (3.8) and (3.9)
that
\begin{eqnarray}
\sum_{|b_{\mu}|<\alpha^{\frac13} r}\log\frac{(\alpha^{\frac13}+1)r}{|z-b_{\mu}|}&\leq&\sum_{\mu=1}^p\log\frac{(\alpha^{\frac13}+1)r}{|z-b_{\mu}|}=\sum_{\mu=1}^p\log\frac{(\alpha^{\frac13}+1)r}{|z-\hat{b_{\mu}}|}\nonumber\\
&\leq&\sum_{\mu=1}^p\log\frac{(\alpha^{\frac13}+1)rp}{B\mu}\leq\sum_{\mu=1}^p\log\frac{\alpha^{\frac13}(\alpha^{\frac13}+1)p}{\varepsilon
\mu}.
\end{eqnarray}
Denote by
$A_2=\frac{\alpha^{\frac13}(\alpha^{\frac13}+1)}{\varepsilon }$ for
simplicity, then by the Stirling's formula and (3.8) we have
\begin{eqnarray}
\sum_{|b_{\mu}|<\alpha^{\frac13} r}\log\frac{(\alpha^{\frac13}+1)r}{|z-b_{\mu}|}
& \leq & \sum_{\mu=1}^p\log\frac{A_2p}{\mu}=
\log \frac{(A_2p)^p}{\Gamma(p+1)}\nonumber\\
&\leq & \log \frac{(A_2p)^p}{A_3p^{\frac12}(\frac{p}{e})^p}=\log \frac{(A_2e)^p}{A_3p^{\frac12}}\nonumber\\
&\leq&\log A_4^p=p\log A_4=\log A_4\cdot n(\alpha^{\frac{h+2}{3}},\phi)
\nonumber\\
&\leq& A_5 n(\alpha^{\frac23}r, \phi),
\end{eqnarray}
here $A_3>0,\ A_4>1,\ A_5>0$ are constants independent of $r$. For $\alpha^{\frac23}r>1$, we have
\begin{eqnarray}
&&N(\alpha r, \phi)\geq \int_{\alpha^{\frac23}r}^{\alpha r}\frac{n(t,\phi)-n(0, \phi)}{t}dt+n(0,\phi)\log(\alpha r)\nonumber\\
&\geq &n(\alpha^{\frac23}r,\phi)\int_{\alpha^{\frac23}r}^{\alpha r}\frac{dt}{t}-n(0,\phi)\int_{\alpha^{\frac23}r}^{\alpha r}\frac{dt}{t}+n(0,\phi)\log(\alpha r)\nonumber\\
&\geq &n(\alpha^{\frac23}r,\phi)\frac{\alpha
r-\alpha^{\frac23}r}{\alpha
r}=(1-\alpha^{-\frac13})n(\alpha^{\frac23}r,\phi),\end{eqnarray}
then \be n(\alpha^{\frac23} r,\phi) \leq
\frac{1}{1-\alpha^{-\frac13}}N(\alpha r,\phi).\ee (3.11) and
(3.13) yields \be \sum_{|b_{\mu}|<\alpha^{\frac13}
r}\log\frac{(\alpha^{\frac13}+1)r}{|z-b_{\mu}|} \leq A_6 N(\alpha
r,\phi)\leq A_6 T(\alpha r,\phi)\ee with constant $A_6>0$.

For each $h$, we define (it has been mentioned that $q$ and $D_1,
D_2, \cdots, D_q$ depend on $h$) \be Y_h=\{ r:\ \mbox{there exist} \
z\in \cup_{j=1}^qD_j\ \mbox{such that} \ |z|=r\},\ee \be E_h=Y_h
\cap [\alpha^{\frac{h}{3}},\alpha^{\frac{h+1}{3}}].\ee Then \be
\int_{E_h}1dx\leq\int_{Y_h}1dx\leq 4B=4\varepsilon
\alpha^{\frac{h}{3}}.\ee Set \be
E=\cup_{h=0}^{\infty}E_h\cap(1,\infty).\ee Then by (3,7) and (3.14),
we have for all $z$ satisfying $|z|=r\notin [0,1]\cup E$, \be
|\phi(z)|\leq e^{AT(\alpha r,\phi)} \ee with $A=A_1+A_6$.

For any $x>1$, there exist nonnegative integer $h$ such that \be
\alpha^{\frac{h}{3}}<x\leq \alpha^{\frac{h+1}{3}}.\ee It follows
form (3.20) and (3.17) that
\begin{eqnarray}  \int_{E\cap (1,x]}\frac{1}{t}dt&\leq& \int_{E\cap (1,\alpha^{\frac{h+1}{3}}]}\frac{1}{t}dt=\sum_{j=0}^h\int_{E_j}\frac{1}{t}dt\leq \sum_{j=0}^h\frac{1}{\alpha^{\frac{j}{3}}}4\varepsilon \alpha^{\frac{j}{3}}\nonumber\\
&=&4\varepsilon (h+1)
\leq 12\varepsilon\frac{\log x}{\log\alpha}+4\varepsilon.\end{eqnarray}
Therefore
\be \delta(E)=\limsup_{x\rightarrow+\infty}\frac{\int_{E\cap (1,x]}\frac{1}{t}dt}{\log x}<\frac{12\varepsilon}{\log\alpha}. \ee
Since $\varepsilon$ is arbitrary small, the proof is completed.\qed

\vskip 3mm

{\bf Lemma 4.} \ Let $N$ be a positive integer. Suppose that the
Dirichlet series \be F_i(s)=\sum_{n=1}^{\infty}\frac{a_i(n)}{n^s},\
\ \ i=1,\ 2,\ \cdots, \ N\ee are convergent in the region
$\sigma_0<\sigma<\infty$, and for each $i=1,\ 2,\ \cdots, \ N$,
$\phi_i(s)$ is meromorphic function in the complex plane with
Nevanlinna characteristic satisfying \be T(r, \phi_i)=o(r) \ \
\mbox{as}\ \ r\rightarrow \infty.\ee Suppose that \be
\sum_{i=1}^N\phi_i(s)F_i(s)=0 \ee holds identically in
$\sigma_0<\sigma<\infty$. Then for each positive integer $n$, \be
\sum_{i=1}^Na_i(n)\phi_i(s)=0 \ee holds identically in the complex
plane.

\vskip 3mm

{\bf Proof.} \  We prove by contradiction. Let $n_0$ be the
minimal index for which \be \sum_{i=1}^Na_i(n_0)\phi_i(s)\ee
 is not identically zero in the complex plane. Then by (3.23) and (3.25), we have
\be
[\sum_{i=1}^Na_i(n_0)\phi_i(s)]\frac{1}{n_0^s}+\sum_{i=1}^N\phi_i(s)[\sum_{n=n_0+1}^{\infty}\frac{a_i(n)}{n^s}]=0
\ee identically in $\sigma_0<\sigma<\infty$. Denote by \be
\psi_i(s)=\frac{\phi_i(s)}{\sum_{i=1}^Na_i(n_0)\phi_i(s)},\ \ \
i=1,\ 2,\ \cdots, \ N. \ee Then we get \be
-1=\sum_{i=1}^N\psi_i(s)n_0^s[\sum_{n=n_0+1}^{\infty}\frac{a_i(n)}{n^s}]
\ee identically in $\sigma_0<\sigma<\infty$. By using the well-known
properties of the Nevanlinna characteristic, we have from (3.24) and
(3.29) that \be T(r, \psi_i)=o(r) \ \ \mbox{as}\ \ r\rightarrow
\infty, \ \ \ i=1,\ 2,\ \cdots, \ N.\ee Then by using Lemma 3 with
$\varepsilon$ small enough, we deduce that there exists real
sequence $\max\{0, \sigma_0\}<\sigma_k\rightarrow\infty$, such that
for each $\epsilon>0$, there exists $k_{\epsilon}>0$ such that for
$k>k_{\epsilon}$, \be  |\psi_i(\sigma_k)|\leq e^{\epsilon \sigma_k},
\ \ \ i=1,\ 2,\ \cdots, \ N.\ee By the general theory of Dirichlet
series, there exist $M>0$ such that \be
|\sum_{n=n_0+1}^{\infty}\frac{a_i(n)}{\sigma_k^s}|\leq
M(n_0+1)^{-\sigma_k},\ \ \ \ k>k_{\epsilon}, \ \ i=1,\ 2,\ \cdots, \
N .\ee By (3.30) (3.32) and (3.33) we have for $k>k_{\epsilon}$, \be
1=|\sum_{i=1}^N\psi_i(\sigma_k)n_0^{\sigma_k}[\sum_{n=n_0+1}^{\infty}\frac{a_i(n)}{n^{\sigma_k}}]|\leq
NMe^{\epsilon
\sigma_k}(\frac{n_0}{n_0+1})^{\sigma_k}=NMe^{(\epsilon-\log\frac{n_0+1}{n_0})\sigma_k}.
\ee Thus we may choose $\epsilon$ such that
$0<\epsilon<\log\frac{n_0+1}{n_0}$, then the right hand side of
(5.34) tends to zero as $k\rightarrow\infty$ and we have a
contradiction. \qed

\vskip 3mm

{\bf Proof of Theorem 1.} \ Assume that $\zeta(s)$ satisfy the
following difference equation \be
\sum_{i=1}^N\phi_i(s)\zeta(s+s_0)^{k_0(i)}\zeta(s+s_1)^{k_1(i)}\cdots\zeta(s+s_m)^{k_m(i)}
=0,\ee where $\phi_i(s), \ i=1,\ 2,\ \cdots, \ N$, are meromorphic
functions not all are identically zero in the complex plane, with
Nevanlinna characteristic satisfying \be T(r, \phi_i)=o(r) \ \
\mbox{as}\ \ r\rightarrow \infty,\ee and $K(i)=(k_0(i), k_1(i),
\cdots, k_m(i)), \ i=1,\ 2,\ \cdots, \ N$, are multi-indices with
every index being a non-negative integer such that \be
K(i)\not=K(j), \ \ \ \ 1\leq i<j\leq N. \ee Here $K(i)=K(j)$ means
that $\ k_l(i)=k_l(j)$ for each $l,\ 0\leq l\leq m$. Let \be
\sum_{n=1}^{\infty}\frac{A_i(n)}{n^s}:=\zeta(s+s_0)^{k_0(i)}\zeta(s+s_1)^{k_1(i)}\cdots\zeta(s+s_m)^{k_m(i)},\
\ \ i=1,\ 2,\ \cdots, \ N,\ee where the Dirichlet series is
convergent in the region $1-\min\{\Re s_0, \Re s_1, \cdots, \Re
s_m\}<\Re s<\infty$. By Lemma 4, we have for each positive integer
$n$, \be \sum_{i=1}^NA_i(n)\phi_i(s)=0 \ee identically in the
complex plane.

Let $S$ be a fixed point such that not all of $\phi_i(S), \ i=1,\
2,\ \cdots, \ N$, are equal to zero, then \be
\sum_{i=1}^NA_i(n)\phi_i(S)=0,\ \ \ \ n=1,\ 2,\ 3, \cdots. \ee We
deduce from (3.40) and (3.38) that \be
\sum_{i=1}^N\phi_i(S)\zeta(s+s_0)^{k_0(i)}\zeta(s+s_1)^{k_1(i)}\cdots\zeta(s+s_m)^{k_m(i)}
=0 \ee holds identically in the complex plane. That is, $\zeta(s)$
satisfies a nontrival difference equation with constant
coefficients. This is a contradiction to Theorem 2. This completes
the proof. \qed

\vskip 8mm

{\bf Acknowledgements.} \ The authors would like to express their
sincere thank to the referee for his valuable comments to our
paper, and in particular for the shortening of the original proof
to Lemma 1., which greatly improved the readability of the paper.

\vskip 8mm

\noindent {\large\bf References} \vskip 8mm

\re{[1]} S. B. Bank and R. P. Kaufman, An extension of
H\"{o}lder's theorem concerning the Gamma function, Funkcial.
Ekvac. {\bf 19}(1976), 53-63.

\re{[2]} H. Cartan, Sur les syst\`{e}mes de fonctions holomorphes
\`{a} vari\'{e}t\'{e}s lin\'{e}aires lacunaires et leurs
applications, Ann. Sci. Ecole Norm. Sup. (3) {\bf 45}(1928),
255-346.

\re{[3]} F. Hausdorff, Zum H\"{o}lderschen Satz \"{u}ber
$\Gamma(z)$, Math. Ann. {\bf 94}(1925), 244-247.

\re{[4]} W. K. Hayman, Meromorphic functions, Oxford Univ. Press,
Oxford, 1964.

\re{[5]} D. Hilbert, Mathematische Probleme, in: {\sl Die
Hilbertschen Probleme}, Leipzig: Akademische Verlagsgesellschaft
Geest \& Portig, 23-80. 1971.

\re{[6]} O. H\"{o}lder, \"{U}ber die Eigenschaft der
$\Gamma$-Function, keiner algebraischen Differentialgleichung zu
gen\"{u}gen, Math. Ann. {\bf 28}(1887), 1-13.

\re{[7]} A. A. Karatsuba and S. M. Voronin, The Riemann
Zeta-Function. Translated from the Russian by N. Koblitz, Walter
de Gruyter, Berlin, New York, 1992.

\re{[8]} B. J. Levin, Distribution of zeros of entire functions,
revised edition, translated from the Russian by R. P. Boas et al.,
American Mathematical Society, Providence, 1980.

\re{[9]} E. Moore, Concerning transcendentally transcendental
functions. Math. Ann. {\bf 48}(1897), 49-74.

\re{[10]} D. D. Mordykhai-Boltovskoi, On hypertranscendence of the
function $\xi(x, s)$, Izv. Politekh. Inst. Warsaw. {\bf 2}(1914),
1-16.

\re{[11]} A. Ostrowski, Neuer Beweis des H\"{o}lderschen Satzes,
dab die $\Gamma$-Function keiner algebraischen
Differentialgleichung gen\"{u}gt. Math. Ann. {\bf 79}(1919),
286-288.

\re{[12]} A. Ostrowski, Zum H\"{o}lderschen Satz \"{u}ber
$\Gamma(x)$, Math. Ann. {\bf 94}(1925), 248-251.

\re{[13]} A. Ostrowski, \"{U}ber Dirichletsche Reihen und
algebraische Differentialgleichungen, Math. Z. {\bf 8}(1920),
241-298.

\re{[14]} A. G. Postnikov, On the differential independence of
Dirichlet series, Dokl. Akad. Nauk SSSR. {\bf 66}(1949), 561-564.

\re{[15]} A. G. Postnikov, Generalization of one of the Hilbert
problems, Dokl. Akad. Nauk SSSR. {\bf 107}(1956), 512-515.

\re{[16]} E. C. Titchmarsh, The Theory of the Riemann Zeta
Function. 2nd edition. Revised by D. R. Heath-Brown, Clarendon
Press, Oxford, 1986.

\re{[17]} S. M. Voronin, The distribution of the nonzero values of
the Riemann zeta-function, Trudy Mat. Inst. Steklov. {\bf
128}(1972), 131-150. 260.

\end{document}